\documentclass[12pt,a4paper]{amsart}
\usepackage{amssymb}
\usepackage{mathrsfs}

\headsep=1cm \numberwithin{equation}{section}
\newtheorem{thm}{Theorem}[section]
\newtheorem{cor}[thm]{Corollary}

\theoremstyle{definition}
\newtheorem{defn}[thm]{Definition}
\theoremstyle{remark}
\newtheorem{rem}[thm]{Remark}
\numberwithin{equation}{section}

\newcommand{\nat}{\mathbb N}
\newcommand{\inte}{\mathbb Z}
\newcommand\Supp{\operatorname{Supp}}
\newcommand\Ass{\operatorname{Ass}}

\newcommand\Spec{\operatorname{Spec}}

\newcommand\Hom{\operatorname{Hom}}
\newcommand\Ext{\operatorname{Ext}}

\newcommand{\p}{\frak p }
\newcommand{\q}{\frak q }

\begin{document} \title [Faltings' local-global principle  for the finiteness of local cohomology]
{Faltings' local-global principle  for the finiteness of local cohomology modules}
\author{Davood Asadollahi and Reza Naghipour$^*$}
\address{Department of Mathematics, University of Tabriz, Tabriz, Iran;
and School of Mathematics, Institute for studies in Theoretical
Phisics and Mathematics (IPM), P.O. Box 19395-5746, Tehran, Iran.}
\email{naghipour@ipm.ir} \email {naghipour@tabrizu.ac.ir}
\email{d\_asadollahi@tabrizu.ac.ir}

\thanks{ 
This research was in part supported by a grant from IPM. \\
$^*$Corresponding author: e-mail: {\it naghipour@ipm.ir} (Reza
Naghipour)}%

\begin{abstract}
Let $(R, \frak m)$ be a complete local ring, $\frak a$ an ideal of $R$ and $M$ a finitely generated $R$-module.
The aim of this paper is to show that for any non-negative integer $n$,
$f^n_{\frak a}(M)={\rm inf} \{0\leq i\in\inte|\, \dim H^{i}_{\frak a}(M)/N \geq n  \text{ for any finitely generated submodule}\,\, N \subseteq  H^{i}_{\frak a}(M)\}$,
where $f^n_{\frak a}(M):=\inf\{f_{\frak a R_{\frak p}}(M_{\frak p})\,\,|\,\,{\frak p}\in \Supp M/\frak a M\,\,{\rm and}\,\,\dim R/{\frak p}\geq n\}$
is  the $n$-th finiteness dimension of $M$ relative to $\frak a$. As a consequence, it follows that the set
$$ \Ass_R(\oplus _{i=0}^{f^n_{\frak a}(M)}H^{i}_{\frak a}(M))\cap \{\frak p\in \Spec R|\, \dim R/\frak p\geq n\}$$ is finite. This generalizes the main result
of Quy \cite{Qu} and  Brodmann-Lashgari \cite{BL}.

\end{abstract}
\maketitle
\section{Introduction}
Throughout this paper, let $R$ denote a commutative Noetherian ring
(with identity) and $\frak a$ an ideal of $R$. For an $R$-module $M$, the
$i^{\rm th}$ local cohomology module of $M$ with support in $V(\frak a)$
is defined as:
$$H^i_{\frak a}(M) = \underset{n\geq1} {\varinjlim}\,\, \Ext^i_R(R/\frak a^n, M).$$  We refer the reader to \cite{BS} or \cite{Gr1} for more
details about local cohomology.  An important theorem in local
cohomology is Faltings' Local-global Principle for the Finiteness
Dimension of local cohomology modules \cite[Satz 1]{Fa1}, which
states that for a positive integer $r$,  the $R_{\frak p}$-module $H^i_{\frak a R_{\frak p}}(M_{\frak p})$ is finitely generated
for all $i\leq r$ and for all ${\frak p}\in  \Spec R$ if and only if
the $R$-module $H^i_{\frak a}(M)$ is finitely generated for all $i\leq r$.

Another formulation of Faltings' Local-global Principle,
particularly relevant for this paper, is in terms of the
generalization of the finiteness dimension $f_{\frak a}(M)$ of $M$ relative
to $\frak a $, where
$$f_{\frak a}(M):=\inf\{i\in \Bbb{N}_0\,\,|\,\,H^i_{\frak a}(M)\,\,{\rm is}\,\,{\rm not}\,\,{\rm finitely}\,\,{\rm generated}
\},\,\,\,\,\,\,\,\,\,\,\,(\dag)$$ with the usual convention that the
infimum of the empty set of integers is interpreted as $\infty$. It
is well known that
$f_{\frak a}(M)=\inf\{f_{\frak a R_{\frak p}}(M_{\frak p})\,\,|\,\,{\frak p}\in \Supp M/\frak aM\,\,{\rm and} \,\,\dim R/{\frak p}\geq 0\}$,
see \cite[9.6.2]{BS}.  Using this idea, for any non-negative integer $n$,   K. Bahmanpour et al., in \cite{BNS1},  introduced the notion of the {\it $n$-th finiteness dimension} $f^n_{\frak a}(M)$ of $M$ relative to $\frak a$ by
 $$f^n_{\frak a}(M):=\inf\{f_{\frak a R_{\frak p}}(M_{\frak p})\,\,|\,\,{\frak p}\in \Supp M/\frak a M\,\,{\rm and}\,\, \dim R/ \frak p\geq n\}.$$
Note that $f^n_{\frak a}(M)$ is either a positive integer or $\infty$ and that $f^0_{\frak a}(M)=f_{\frak a}(M)$.  Also, they showed that the least integer $i$ such that $H^i_{\frak a} (M)$ is not minimax (resp.  weakly Laskerian), equals to $f^1_{\frak a} (M)$ (resp. $f^2_{\frak a} (M)$). So it is rather
natural to ask whether Faltings' Local-global Principle, as stated in $(\dag)$, generalizes in the obvious way to the invariants  $f^n_{\frak a}(M)$.
To this end, we are able to introduce the concept of an {\it $R$-module in dimension $< n$} as a generalization of the notion of a FSF module  \cite{Qu}.
An $R$-module $M$ is said to be in dimension $< n$, if there is a finitely generated submodule $N$ of $M$ such that $\dim \Supp  M/ N <n$.

As a main result of this paper we shall show that:

\begin{thm}
Let $(R,{\frak m})$ be a complete local ring, ${\frak a}$  an ideal of $R$ and $M$ a finitely generated $R$-module. Then for any $n\in\mathbb{N}_0$,
$$f_{\frak a}^n(M)={\rm inf} \{0\leq i\in\inte|\, H^{i}_{\frak a}(M) \text{is not in dimension}\, < n\}. $$
\end{thm}

As a consequence of Theorem 1.2, we derive the following, which is a generalization of the main result of Quy \cite[Theorem 3.2]{Qu}  and  Brodmann-Lashgari \cite[Theorem 2.2]{BL}.

\begin{cor}
Let $(R,{\frak m})$ be a complete local ring, ${\frak a}$  an ideal of $R$ and $M$ a finitely generated $R$-module. Then
the set
\begin{center}
$ \Ass_R(\bigoplus _{i=0}^{f^n_{\frak a}(M)}H^{i}_{\frak a}(M))\cap \{\frak p\in \Spec R|\, \dim R/\frak p\geq n\}$
\end{center}
is finite.
\end{cor}

Throughout this paper, $R$ will always be a commutative Noetherian
ring with non-zero identity and $\frak a$ will be an ideal of $R$.  By a {\it skinny or weakly
Laskerian} module, we mean an $R$-module $M$ such that the set $\Ass_R M/N$ is finite, for each submodule $N$ of $M$  (cf. \cite{Ro} or \cite{DM1}). Moreover, an $R$-module $M$ is said to be {\it minimax}, if there exists a finitely generated submodule $N$ of $M$, such that $M/N$ is Artinian. \\

\section{The Main  Results}
In \cite{Qu}, P. H. Quy introduced the class of FSF modules and he has given some properties and applications of this modules. An $R$-module $M$ is said to be a FSF module if there is a finitely generated submodule $N$ of $M$ such that support of the quotient module $M/N$ is finite. When $R$ is a Noetherian ring, it is clear that, if $M$ is FSF, then $\dim\Supp M/N\leq 1$. This motivates the definition.

\begin{defn}
Let $n$ be a non-negative integer. An $R$-module $M$ is said to be in dimension $< n$, if there is a finitely generated submodule $N$ of $M$ such that $\dim \Supp  M/ N <n$
\end{defn}
\begin{rem}\label{rem}
Let $n$ be a non-negative integer and let $M$ be an $R$-module.

{\rm (i)}  If $n=0$, then $M$ is in dimension $< n$ if and only if $M$ is Noetherian.

{\rm (ii)} If $M$ is minimax, then $M$ is in dimension $< 1$. In particular, if $M$ is Noetherian or Artinian, then $M$ is in dimension $< 1$.

{\rm (iii)}  If $M$ is FSF, then $M$ is in dimension $< 2$.

{\rm (iv)} If $M$ is skinny or weakly Laskerian, then $M$ is in dimension $< 2$, by \cite [Theorem 3.3]{Ba}.

{\rm (v)} If $M$ is reflexive, then $M$ is in dimension $<1$.

{\rm (vi)} If $M$ is linearly compact, then $M$ is in dimension $<1$.  Recall that $M$ is said to be linearly compact if
each system of congruences $x \equiv x_i(M_i)$ indexed by a set $I$, and where the $M_i$ are submodules of $M$,
has a solution $x$ whenever it has a solution for every finite subsystem. It is know that it either $M$ reflexive or  linearly compact,
then $M$ is minimax (see e.g. \cite{DM1}).
\end{rem}

\begin{defn}
If $T$ is a subset of $\Spec R$ and $n\in\nat_0$, then we define $$(T)_{\geq n}=\{ \p\in T | \dim R/\p\geq n\}.$$
\end{defn}
\begin{defn}
Let $n$ be a non-negative integer, $\frak a$  an ideal of $R$ and $M$ an $R$-module. Then we define
$$h_{\frak a}^n(M)={\rm inf} \{0\leq i\in\inte: \,\, H^{i}_{\frak a}(M) \,\, \text{is not in dimension}< n\}. $$
\end{defn}

Now we are prepared to state and prove the  main result of this paper, which shows that the least integer $i$ such that $H^i_{\frak a}(M)$ is not
in dimension $<n$, equals to $\inf \{f_{\frak aR_{\frak p}}(M_{\frak p})\,\,|\,\,\frak p\in \Supp M/\frak aM\,\,{\rm and} \,\,\dim R/\frak
p\geq n\}.$

\begin{thm}
Let $(R,{\frak m})$ be a complete local ring, ${\frak a}$  an ideal of $R$ and $M$ a finitely generated $R$-module. Then for any $n\in\mathbb{N}_0$,
$$f_{\frak a}^n(M)=h_{\frak a}^n(M).$$
\end{thm}

\proof
 Let $i$ be a non-negative integer such that  $ H^{i}_{\frak a}(M)$ is in dimension $< n$. Then it follows from the definition that there is a finitely generated submodule $N$ of $ H^{i}_{\frak a}(M)$ such that
$\dim\Supp H^{i}_{\frak a}(M)/N <n$. Thus for all $\frak p\in\Supp M/{\frak a}M$ with $\dim R/\frak p\geq n$ we have $(H_{\frak a}^i(M)/N)_{\frak p}=0 $. Therefore $(H^{i}_{\frak a}(M))_{\frak p}$ is finitely generated and so $h_{\frak a}^n(M)\leq f_{\frak a}^n(M)$.

We now suppose that $t=h_{\frak a}^n(M)< f_{\frak a}^n(M)$, and look for a contradiction.  To this end, first we show that for all $m\in\mathbb{N}$, the set
$(\Ass_R H^{t}_{\frak a}(M)/{(0:_{H^{t}_{\frak a}(M)}{\frak a}^m)})_{\geq n}$ is finite. To achieve this, suppose the contrary is true. Then there exists a countably infinite subset $\{ \frak p_k\}^{\infty}_{k=1}$ of $(\Ass_R H^{t}_{\frak a}(M)/(0:_{ H^{t}_{\frak a}(M)}{\frak a}^m))_{\geq n}$. Let $S$ be the multiplicatively closed subset $R\backslash \bigcup^{\infty}_{k=1}\p_k$. We now show that the $S^{-1}R$-module $H^{t}_{S^{-1}{\frak a}}(S^{-1}M)$ is finitely generated.  To do this, in view of  Faltings' Local-global Principle theorem (see \cite[Theorem 9.6.1]{BS}), it is enough to show that for all $j\leq t$ and for all prime ideals $\frak p$ with $S\cap\frak p=\emptyset$, the $R_{\frak p}$-module $(H^{j}_{S^{-1}{\frak a}}(S^{-1}M))_{S^{-1}\frak p}$ is finitely generated. Since $S\cap\frak p=\emptyset$, it follows that $\frak p\subseteq \bigcup^{\infty}_{k=1}\p_k$, and  so  by \cite[Lemma 3.2]{MV}
there exists $k\geq1$ such that $\frak p\subseteq \frak p_k$. Thus $\dim R/\frak p\geq n$. Now, as $j< f_{\frak a}^n(M)$,  it follows that
 $(H^{j}_{S^{-1}{\frak a}}(S^{-1}M))_{S^{-1}\p}\cong (H^{j}_{\frak a}(M))_{\p}$ is finitely generated, as required. Therefore
the set  $\Ass_{S^{-1}R} S^{-1} (H^{t}_{\frak a}(M)/{(0:_{H^{t}_{\frak a}(M)}\frak a}^m))$ is finite.
On the other hand, we have  $S^{-1}\p_k\in \Ass_{S^{-1}R}(S^{-1} H^{t}_{\frak a}(M)/{(0:_{H^{t}_{\frak a}(M)}\frak a^m)})$ for all $k=1,2, \dots$, which is a contradiction.

Consequently, for all $m\in\mathbb{N}$, the set $(\Ass_R  H^{t}_{\frak a}(M)/{(0:_{H^{t}_{\frak a}(M)}{\frak a}^m)})_{\geq n}$ is finite. Now, we let
 $\mathbb{A}$ be the set of all prime ideals $\p$ of $R$ such that there exists $m\in\mathbb{N}$  with
 $\p\in(\Ass_R H^{t}_{\frak a}(M)/(0:_{H^{t}_{\frak a}(M)}{\frak a}^m))_{\geq n}$. Then $\mathbb{A}$ is a countably infinite set. Let $S$ be the multiplicatively closed subset $R\backslash \bigcup_{\p\in \mathbb{A}}\p$. Then it is easy to see that  for all $m\in\mathbb{N}$, the set $\Ass_{S^{-1}R}S^{-1} (H^{t}_{\frak a}(M)/(0:_{H^{t}_{\frak a}(M)}{\frak a}^m))$ is finite. Thus for all $m\in\mathbb{N}$, the set 
 $\Supp S^{-1}(H^{t}_{\frak a}(M)/{(0:_{H^{t}_{\frak a}(M)}{\frak a}^m)})$ is a closed subset of $\Spec R$ (in the  Zariski topology), and so the descending chain
 $$\cdots\supseteq \Supp S^{-1}(H^{t}_{\frak a}(M)/(0:_{ H^{t}_{\frak a}(M)}\frak a^m))\supseteq \Supp S^{-1} (H^{t}_{\frak a}(M)/(0:_{ H^{t}_{\frak a}(M)}\frak a^{m+1}))\supseteq \cdots, $$
 is eventually stationary. Let $E_m$ denote its eventually stationary value, so that there is $m\in\mathbb{N}$ such that
 for all $l\geq m$, $$E_m=\Supp S^{-1}(H^{t}_{\frak a}(M)/(0:_{H^{t}_{\frak a}(M)}\frak a^l)).$$

 Since for all $i<t$, $H^{i}_{\frak a}(M)$ is in dimension $<n$, it follows that  $(0:_{H^{t}_{\frak a}(M)}\frak a ^m)$ is also in dimension $<n$. Thus there is a finitely generated submodule $N$ of $(0:_{H^{t}_{\frak a}(M)}{\frak a}^m)$ such that $\dim\Supp\, (0:_{H^{t}_{\frak a}(M)}{\frak a}^m)/N<n$.

Now,  we show that $\dim\Supp H^{t}_{\frak a}(M)/(0:_{H^{t}_{\frak a}(M)}\frak a^m)<n$. Suppose the contrary is true. Then there exists $\q\in\Ass_RH^{t}_{\frak a}(M)/(0:_{H^{t}_{\frak a}(M)}\frak a^m)$ such that $\dim R/\q\geq n$. Whence $\frak q\in \mathbb{A}$, and so
$S\cap\frak q=\emptyset$.  Thus $S^{-1}\q\in E_m$.  On the other hand, since $t<f_{\frak a}^n(M)$ and $\dim R/\q\geq n$, it yields that
$(H^{t}_{\frak a}(M))_{\q}$ is a finitely generated $R_{\frak q}$-module. Therefore there exists $l\geq m$
such that $(\frak aR_{\q})^l(H^{t}_{\frak a}(M))_{\q}=0$, and so $(H^{t}_{\frak a}(M)/(0:_{ H^{t}_{\frak a}(M)}{\frak a}^l))_{\q}=0$. Hence we have $(S^{-1}(H^{t}_{\frak a}(M)/(0:_{H^{t}_{\frak a}(M)}{\frak a}^l)))_{S^{-1}\q}=0$, and so $S^{-1}\q\not \in E_m$, which is a contradiction.
Consequently, we have $\dim\Supp H^{t}_{\frak a}(M)/(0:_{H^{t}_{\frak a}(M)}\frak a^m)<n$. Finally, from the exact sequence
 $$0\longrightarrow (0:_{H^{t}_{\frak a}(M)}\frak a^m)/N \longrightarrow H^{t}_{\frak a}(M)/N\longrightarrow H^{t}_{\frak a}(M)/(0:_{H^{t}_{\frak a}(M)}{\frak a}^m)\longrightarrow 0,$$
we conclude that $\dim\Supp H^{t}_{\frak a}(M)/N<n$. That is the $R$-module $H^{t}_{\frak a}(M)$ is in dimension $< n$, and so we have obtained  a contradiction.  \qed \\

\begin{cor}
Let $(R, \frak m)$ be a complete  local ring, $\frak a$ an ideal of $R$ and $M$ a finitely generated $R$-module. Then

{\rm (i)}  $f^1_I(M)=\inf\{i\in \Bbb{N}_0\,\,|\,\,H^i_I(M)\,\,{\rm is}\,\,{\rm not}\,\,{\rm minimax}\}.$

{\rm (ii)} $f^2_I(M)=\inf\{i\in \Bbb{N}_0\,\,|\,\,H^i_I(M)\,\,{\rm is}\,\,{\rm not}\,\,{\rm weakly}\,\,{\rm Laskerian}\}.$
\end{cor}

\proof The result follows from the definition and Theorem 2.5. \qed\\

\begin{cor}
Let $(R, \frak m)$ be a complete  local ring, $\frak a$ an
ideal of $R$ and $M$ a finitely generated $R$-module.
Then the set $(\Ass_R(\oplus_{i=0}^{t-1}H^i_{\frak a}(M)))_{\geq n}$ is finite, where $t=f^n_{\frak a}(M)$.
\end{cor}
\proof  In view of Theorem 2.5, the $R$-module $\oplus_{i=0}^{t-1}H^i_{\frak a}(M)$ is in dimension $<n$. So there is a finitely generated submodule $N$ of $\oplus_{i=0}^{t-1}H^i_{\frak a}(M)$ such that $\dim \oplus_{i=0}^{t-1}H^i_{\frak a}(M)/N\leq n-1$. Now the assertion follows immediately from
 the exact sequence $$0 \longrightarrow N \longrightarrow \oplus_{i=0}^{t-1}H^i_{\frak a}(M) \longrightarrow \oplus_{i=0}^{t-1}H^i_{\frak a}(M)/N \longrightarrow 0.$$  \qed

The final result is a generalization of the main result of Quy \cite[Theorem 3.2]{Qu}  and  Brodmann-Lashgari  \cite[Theorem 2.2]{BL} for complete local rings.
 \begin{thm}
Let $(R,{\frak m})$ be a complete local ring, $\frak a$ an
ideal of $R$ and $M$ a finitely generated $R$-module.
Then the set $(\Ass_R(\oplus_{i=0}^{t}H^i_{\frak a}(M)))_{\geq n}$ is  finite, where $t=f^n_{\frak a}(M)$.
\end{thm}

\proof Using Corollary 2.7 it is enough to show that the set $(\Ass_RH^t_{\frak a}(M))_{\geq n}$ is finite. To this end, suppose that the contrary is true. Then  there exists a countably infinite subset
$\{\frak q_{k}\}_{k=1}^{\infty}$ of $(\Ass_RH^t_{\frak a}(M))_{\geq n}$, and so  by \cite[Lemma 3.2]{MV},
$\mathfrak{p} \not\subseteq \bigcup_{k=1}^{\infty}\frak q_{k}$, for every $\frak p\in \Spec R$ with $\dim R/\frak p<n$. Let
$S$ be the multiplicatively closed subset $R \backslash \bigcup
_{k=1}^{\infty}\frak q_{k}$. Then, it easily follows from
\cite[Theorem 9.6.1]{BS} that $S^{-1}H^i_{\frak a}(M)$ is a finitely generated
$S^{-1}R$-module, for all $i=0,1,\dots,t-1$ and so in view of \cite[Theorem 2.3]{BN}, the $R$-module $$\Hom_{S^{-1}R}(S^{-1}R/S^{-1}\frak a, S^{-1}H^t_{\frak a}(M)),$$ is finitely generated. Therefore the set $$\Ass_{S^{-1}R} S^{-1}H^t_{\frak a}(M)$$ is  finite.
But $ S^{-1}\frak q_{k} \in \Ass_{S^{-1}R} S^{-1}(H^i_{\frak a}(M))$ for all
$k=1,2,\dots, $ which is a contradiction.\qed

\begin{center}
{\bf Acknowledgments}
\end{center}
The authors would like to thank Professors Hossein Zakeri and Kamal Bahmanpour for reading of the first draft and valuable discussions.
Also, we  would like to thank from the Institute for Research in Fundamental Sciences (IPM), for its financial support.

\end{document}